\documentclass[a4paper,11pt,reqno,twoside]{article}
\usepackage{amssymb}
\usepackage[frenchb,english]{babel}
\usepackage{amscd}
\usepackage{amsmath,amsthm,amsfonts,amssymb,graphicx}
\usepackage[colorlinks,linkcolor=blue,anchorcolor=blue,citecolor=green]{hyperref}
\usepackage{mathrsfs}
\usepackage{fancyhdr}
\usepackage{multicol}
\usepackage{abstract}
\usepackage{mathrsfs}
\usepackage[all, pdf]{xy}
\usepackage{bm}

\begin{document}
\newtheorem{lem}{Lemma}[section]
\newtheorem{thm}{Theorem}[section]
\newtheorem{prop}{Proposition}[section]
\newtheorem{cor}{Corollary}[section]
\newtheorem{rem}{Remark}[section]

\noindent  {\Large
$L^{p}$ regularity of weighted Bergman projection on Fock-Bargmann-Hartogs domain}\\

\noindent\text{Le He, Yanyan Tang \; \& \; Zhenhan Tu$^{*}$ }\\

\noindent\small {School of Mathematics and Statistics, Wuhan
University, Wuhan, Hubei 430072, P.R. China}

\noindent\text{Email: hele2014@whu.edu.cn (L. He), yanyantang@whu.edu.cn (Y. Tang), zhhtu.math@whu.edu.cn (Z. Tu)}

\renewcommand{\thefootnote}{{}}
\footnote{\hskip -16pt {$^{*}$Corresponding author. \\}}\\

\normalsize \noindent\textbf{Abstract}\quad
The Fock-Bargmann-Hartogs domain $D_{n, m}(\mu)$ is defined by
$$ D_{n, m}(\mu):=\{(z, w)\in\mathbb{C}^{n}\times\mathbb{C}^m:\Vert w \Vert^2<e^{-\mu\Vert z \Vert^2}\},$$
where $\mu>0.$  The Fock-Bargmann-Hartogs domain $D_{n, m}(\mu)$ is an unbounded strongly pseudoconvex domain with smooth real-analytic boundary.
In this paper, we first compute the weighted Bergman kernel of $D_{n, m}(\mu)$ with respect to the weight $(-\rho)^{\alpha}$, where $\rho(z,w):=\|w\|^2-e^{-\mu \|z\|^2}$ is a defining function for $D_{n, m}(\mu)$ and $\alpha>-1$. Then, for $p\in [1,\infty),$ we show that the
corresponding weighted Bergman projection $P_{D_{n, m}(\mu), (-\rho)^{\alpha}}$ is unbounded on $L^p(D_{n, m}(\mu), (-\rho)^{\alpha})$, except for the trivial case $p=2$.
In particular,  this paper gives an example of an unbounded strongly pseudoconvex domain whose ordinary Bergman projection is $L^p$ irregular when $p\in [1,\infty)\setminus\{2\}$. This result turns out to be completely different from the well-known positive $L^p$ regularity result on bounded strongly pseudoconvex domain.

\vskip 10pt

\noindent \textbf{Keywords:} Fock-Bargmann-Hartogs domain \textperiodcentered \; $L^p$ regularity \textperiodcentered \; Weighted Bergman kernel \textperiodcentered\;Weighted Bergman projection
\vskip 10pt

\noindent \textbf{Mathematics Subject Classification (2010):}  32A36 \textperiodcentered \, 32A25   \textperiodcentered \, 32A07

\pagenumbering{arabic}
\renewcommand{\theequation}
{\arabic{section}.\arabic{equation}}

\setcounter{section}{0}
\setcounter{equation}{0}
\section{Introduction}
\subsection{Setup and Problems}
Let $\Omega$ be a domain in $\mathbb{C}^n$ and $\eta(z)$ be a non-negative measurable function on $\Omega$. For $p\in [1,+\infty)$, let $L^p(\Omega, \eta)$ denote the set of all complex measurable functions $f$ with
$$\bigg(\int_\Omega |f(z)|^p\eta(z)dV(z)\bigg)^{1/p}<+\infty,$$
where $dV(z)$ is the ordinary Lebesgue measure on $\Omega$.
We call $\eta(z)$ a weight on $\Omega$ and $L^p(\Omega, \eta)$ the weighted $L^p$ space of $\Omega$. The norm of $L^p(\Omega, \eta)$ is defined as follows
$$\|f\|_{p,\eta}=\bigg(\int_{\Omega}|f(z)|^p\eta(z) dV(z)\bigg) ^{1/p}.$$
If $p=2$, $L^2(\Omega, \eta)$ is a Hilbert space with the inner product:
$$\left<f,g\right>_\eta=\int_{\Omega}f(z)\overline{g(z)}\eta(z) dV(z).$$
The weighted Bergman space of $\Omega$ with weight $\eta$ is defined by
$$A^p(\Omega, \eta):=\mathcal{O}(\Omega)\cap L^p(\Omega, \eta),$$
where $\mathcal{O}(\Omega)$ is the space of all holomorphic functions on $\Omega$.
Thus $A^2(\Omega, \eta)$ is a subspace of holomorphic functions in $L^2(\Omega, \eta)$. By the result in \cite{Pas}, we know: if $\eta$ is continuous and never vanishes inside $\Omega$, then $A^p(\Omega, \eta)$ is a closed subspace of $L^2(\Omega, \eta)$, and there exists the the orthogonal projection, called the weighted Bergman projection:
$$P_{\Omega, \eta}: L^2(\Omega, \eta)\to A^2(\Omega, \eta).$$
This projection is an integral operator with the weighted Bergman kernel, denoted by $K_{\Omega, \eta}(z,w)$:
$$P_{\Omega, \eta}f(z):=\int_{\Omega}K_{\Omega, \eta}(z,w)f(w)\eta(w) dV(w).$$
When $\eta(z)\equiv 1$, the weighted Bergman kernel $K_{\Omega, \eta}$ and the weighted Bergman projection $P_{\Omega, \eta}$ will degenerate to the ordinary Bergman kernel $K_\Omega$ and the ordinary Bergman projection $P_\Omega$, respectively.

For an arbitrary domain $\Omega\subset\mathbb{C}^n$, a continuous and positive weight $\eta$ on $\Omega$, the corresponding weighted Bergman projection $P_{\Omega, \eta}$ is originally defined on $L^2(\Omega, \eta)$, mapping onto the weighted Bergman space $A^2(\Omega, \eta)$. When we say the weighted Bergman projection $P_{\Omega, \eta}$ on $L^p(\Omega, \eta)$ will mean $P_{\Omega, \eta}$ on the subspace $L^p(\Omega, \eta)\cap L^2(\Omega, \eta)$ of $L^p(\Omega, \eta)$. Therefore, for any $p\in[1,\infty)$, when we say the weighted Bergman projection $P_{\Omega, \eta}$ is bounded on $L^p(\Omega, \eta)$, we mean the weighted Bergman projection $P_{\Omega, \eta}$ mapping $L^p(\Omega, \eta)\cap L^2(\Omega, \eta)$ onto $A^p(\Omega, \eta)\cap L^2(\Omega, \eta)$ is bounded.

Fixing now a domain $\Omega$ and a positive continuous weight $\eta$ on $\Omega$, we define the operator norm of $P_{\Omega, \eta}$ as follows
$$\|P_{\Omega,\eta}\|_{p,\eta}:=\sup\Big\{\frac{\| P_{\Omega,\eta}f\|_{p, \eta}}{\| f\|_{p,\eta}}: f\in L^p(\Omega, \eta)\cap L^2(\Omega, \eta), f\neq 0\Big\}.$$
It is easy to see that $\Vert P_{\Omega, \eta}\Vert_{2,\eta}=1$ in the case of $p=2$. A natural and interesting question is to determine the range of $p\in (1, +\infty)$ such that the weighted Bergman projection $P_{\Omega, \eta}$ is bounded on $L^p(\Omega, \eta)$, except for the trivial case $p=2$.
The question mentioned above is the so-called $L^p$ regularity problem.
\subsection{Background}
The $L^p$ regularity of the (weighted) Bergman projection is a fact of interesting and fundamental importance. Even though two domains are biholomorphic equivalence, the corresponding $L^p$ behavior of the Bergman projection on these two domains may be quite different from each other. There are many papers considering this problem in different settings. One of the most common object is the bounded domain with various boundary conditions. For example, positive $L^p$ regularity results have been obtained on the following domains for all $p\in (1, +\infty)$ in the unweighted version:
\begin{enumerate}
  \item [-] $\Omega$ is bounded strongly pseudoconvex domain (see Lanzani-Stein \cite{Lan2}, Phong-Stein \cite{Pho}).
  \item [-] $\Omega$ is bounded, smooth and pseudoconvex domain of finite type in $\mathbb{C}^{2}$ (see McNeal \cite{Mc1}).
  \item [-] $\Omega$ is bounded, smooth and convex domain of finite type in $\mathbb{C}^{n}$ (see McNeal \cite{Mc2}, McNeal-Stein \cite{Mc3}).
  \end{enumerate}
  We refer to Charpentier-Dupain \cite{Char} and Huo \cite{Huo} for the positive results on other bounded domians. There are examples of smoothly bounded pseudoconvex domains where the $L^p$ boundedness does not hold on the full interval $(1, +\infty)$ (see Barrett-\c{S}ahuto\u{g}lu \cite{Bar}). In addition, if the domain $\Omega$ has serious boundary singularity, in general, there will be a restricted range of p for the $L^p$ boundedness of $P_{\Omega}$ (see, e.g., Chakrabarti-Zeytuncu \cite{Chak}, Edholm-McNeal \cite{Edh}). Many authors are also interested in what happens for the $L^p$ boundedness of weighted Bergman projections if they introduce different type of weights on bounded domains, e.g.,
 \begin{enumerate}
   \item[-] Unit ball $\mathbb{B}^n$ in $\mathbb{C}^n~(n\geq1)$. Let $\omega=(-\rho)^{\alpha}~(\alpha>-1)$, where $\rho(z)=\Vert z \Vert^2-1$ is the defining function of $\mathbb{B}^n$. It is shown that the weighted Bergman projection $P_{\mathbb{B}^n, \omega}$ is bounded from $L^p(\mathbb{B}^n, \omega)$ to $A^p(\mathbb{B}^n, \omega)$ for any $p\in (1, +\infty)$. The result implies that the $L^p$ boundedness is independent on the parameter $\alpha$ (see Rudin \cite[Section 7.1]{Ru}).
   \item[-] Hartogs triangle $\mathbb{H}$ in $\mathbb{C}^2$. Let $\omega(z):=|z_2|^s~(s\in \mathbb{R}), z=(z_1, z_2)\in \mathbb{H}$. Then the range of $p$ for the $L^p$ boundedness of $P_{\mathbb{H}, \omega}$ is related to the power $s$ (see Chen \cite{Che1}).
\end{enumerate}
For more results on other bounded domains with exponential weights, we refer to \v{C}u\v{c}kovi\'{c}-Zeytuncu \cite{Cu2},  Zeytuncu \cite{Zey} and references therein.

For the unbounded domains, very few, however, have studied the $L^p$ regularity problem on this class of domains.
In \cite{Kra2}, Krantz and Peloso determined the $L^p$-mapping properties of the Bergman projection on unbounded, non-smooth worm domians, which is facilitated by the fact that the boundaries of these domains are Levi flat. In \cite{Jan}, Janson, Peetre and Rochberg determined the $L^p$-mapping properties of the Bergman projection on $L^p$ space on $\mathbb{C}^n$ with respect to Gaussian weights $\eta_{\alpha}(z)=e^{-\alpha\|z\|^2}$.
In \cite{Bom}, Bommier-Hato,  Engli\v{s} and Youssfi give criteria for boundedness of the associated Bergman-type projections on $L^p$ space on $\mathbb{C}^n$ with respect to generalized Gaussian weights $e^{-\alpha\|z\|^{2m}}$, where $m>0$.

In this paper, inspired by the above works, we focus on the $L^p$ regularity problem on the Fock-Bargmann-Hartogs domain $D_{n, m}(\mu)$ in $\mathbb{C}^{n+m}$. It is therefore interesting, as a model for the unbounded case, to study the behavior of the (weighted) Bergman projecton on $D_{n, m}(\mu)$.

\subsection{The Fock-Bargmann-Hartogs domain}

For a given positive real number $\mu$, the Fock-Bargmann-Hartogs domain $D_{n,m}(\mu)$ is a Hartogs domain over $\mathbb{C}^n$ defined by
$$D_{n,m}(\mu):=\bigg\{(z,w)\in \mathbb{C}^{n+m}: \|w\|^2<e^{-\mu \|z\|^2}\bigg\},$$
where $\|\cdot\|$ is the standard Hermitian norm. The domain is an unbounded, inhomogeneous strongly pseudoconvex domain in $\mathbb{C}^{n+m}$ with smooth real-analytic boundary. Besides, since each $D_{n, m}(\mu)$ contains $\{(z,0)\in\mathbb{C}^n\times \mathbb{C}^m\} \cong \mathbb{C}^n$, it is also a domain which is not hyperbolic in the sense of Kobayashi. Therefore, it can not be biholomorphic to any bounded domain in $\mathbb{C}^{n+m}$.
For more information of the Fock-Bargmann-Hartogs domain, see Bi-Feng-Tu \cite{Bi}, Kim-Ninh-Yamamori \cite{Kim}, Tu-Wang \cite{Tu} and Yamamori \cite{Yam} and reference therein.

We note that the Fock-Bargmann-Hartogs domain $D_{n, m}(\mu)$ is defined as a domain in $\mathbb{C}^{m+n}$ with the fiber over $\mathbb{C}^n$ being a $m$-dimensional ball. Thus, we could relate the weighted Bergman kernel of $D_{n, m}(\mu)$ to weighted Bergman kernel of the base space $\mathbb{C}^n$ under some condition. Thanks to the relationship between the weighted Bergman kernels, one can get the $L^p$ regular behavior of the corresponding weighted Bergman projections.

\subsection{Main results}
Let
$\rho(z,w)=\|w\|^2-e^{-\mu \|z\|^2}$, $(z,w)\in D_{n,m}(\mu)$. For $-1<\alpha<\infty$, the weighed Bergman space $A^2(D_{n,m}(\mu),(-\rho)^\alpha)$ is defined by
$$A^2(D_{n,m}(\mu),(-\rho)^\alpha):=\bigg\{f\in \mathcal{O}(D_{n,m}(\mu)):\int_{D_{n,m}(\mu)}|f|^2 (-\rho)^{\alpha} dV<\infty\bigg\}.$$
The Bergman kernel of $A^2(D_{n,m}(\mu),(-\rho)^\alpha)$ is denoted by $K_{D_{n,m}(\mu), (-\rho)^\alpha}$. By applying Ligocka's Theorem \cite{Lig}, Yamamori \cite{Yam} gave an explicit expression of the Bergman kernel of $A^2(D_{n,m}(\mu))$.

Firstly, following the method in Bi-Feng-Tu \cite{Bi}, we give a formula of the weighted Bergman kernel $K_{D_{n,m}(\mu), (-\rho)^\alpha}$ as follows.

\begin{thm}\label{thm2}
Let $D_{n,m}(\mu)$ be the Fock-Bargmann-Hartogs domain, $\rho(z,w)=\|w\|^2-e^{-\mu \|z\|^2}$, $(z,w)\in D_{n,m}(\mu)$.
Then, for $\alpha>-1$, the Bergman kernel of the weighted Hilbert space $A^2(D_{n,m}(\mu),(-\rho)^\alpha)$ defined by
\begin{equation*}
A^2(D_{n,m}(\mu),(-\rho)^\alpha):=\bigg\{f\in \mathcal{O}(D_{n,m}(\mu)):\int_{D_{n,m}(\mu)}|f|^2 (-\rho)^{\alpha} dV<\infty\bigg\}
\end{equation*}
can be expressed as
\begin{equation*}
K_{D_{n,m}(\mu), (-\rho)^\alpha}((x,y),(s,t))=\frac{\mu^n}{\pi^{n+m}}\sum_{k\in\mathbb{N}}
\frac{\Gamma(\alpha+m+k+1)(\alpha+m+k)^n}{\Gamma(\alpha+1)\Gamma(k+1)}e^{\mu(\alpha+m+k) \left< x,s \right>}\left< y,t \right>^{k}.
\end{equation*}
\end{thm}

As an application of Theorem \ref{thm2},  we derive a key relation between the weighted Bergman kernel $K_{D_{n, m}(\mu), (-\rho)^\alpha}$ of $D_{n, m}(\mu)$ and the weighted Bergman kernel $K_{\mathbb{C}^n, \eta_{\mu(\alpha+m)}}$ of the base space $\mathbb{C}^n$ (see Lemma \ref{pro:1}).
Next we turn our attention to study the $L^p$ regularity properties of the weighted Bergman projection $P_{D_{n,m}(\mu),(-\rho)^\alpha}$ on $D_{n,m}(\mu)$ with this key relation, and we have the following result.
\begin{thm}\label{th:1}
Let $P_{D_{n,m}(\mu),(-\rho)^\alpha}$ be the weighted Bergman projection on $D_{n,m}(\mu)$ with the weight $(-\rho)^{\alpha}$, where $\rho(z,w)=\|w\|^2-e^{-\mu \|z\|^2}$, $(z,w)\in D_{n,m}(\mu)$ and $\alpha> -1$.
Then, for $1\le p<\infty$, we have that $P_{D_{n,m}(\mu),(-\rho)^\alpha}$ is bounded on $L^p(D_{n, m}(\mu), (-\rho)^\alpha)$ if and only if $p=2$.
\end{thm}

\begin{rem}
Setting $\alpha=0$ in Theorem \ref{th:1}, we obtain that the ordinary Bergman projection $P_{D_{n, m}(\mu)}$ is bounded on $L^p(D_{n, m}(\mu))$ if and only if $p=2$.
\end{rem}

Our proof of Theorems \ref{th:1} employs the technique used by \v{C}u\v{c}kovi\'{c}-Zeytuncu \cite{Cu2} and Zeytuncu \cite{Zey}.
Since $D_{n, m}(\mu)$ is an unbounded domain and it can not be biholomorphic to any bounded domain in $\mathbb{C}^{n+m}$, the Fock-Bargmann-Hartogs domain $D_{n, m}(\mu)$ studied here is different from the bounded Hartogs domain $\Omega_\varphi$ in $\mathbb{C}^2$. What's more, we give an example of an unbounded strongly pseudoconvex domain whose ordinary Bergman projection is $L^p$ irregularuty except for the trivial case $p = 2$.

\section{Preliminaries}
\begin{lem}\label{lem:1}
For $\alpha>-1$, then the following multiple integration exists and
\begin{equation*}
\int_{0}^{1} dx_{m}\cdots \int_{0}^{1-\sum\limits_{i=2}^{m} x_{i}} \bigg(1-\sum\limits_{i=1}^{m}x_{i}\bigg)^{\alpha}\prod\limits_{i=1}^{m}x_{i}^{q_{i}} dx_{1}=\frac{\prod_{i=1}^{m}\Gamma(q_{i}+1)\Gamma(\alpha+1)}{\Gamma(\alpha+\sum_{i=1}^{m}q_{i}+m+1)},
\end{equation*}
where $q=(q_{1},\cdots ,q_{m})\in(\mathbb{R}_{+})^{m}$, here $\mathbb{R}_{+}$ denotes the set of positive real numbers.
\end{lem}

{\bf Proof.}  By  calculating, we have
\begin{equation*}
\begin{aligned}
&\int_{0}^{1} dx_{m}\cdots \int_{0}^{1-\sum\limits_{i=2}^{m} x_{i}} \bigg(1-\sum\limits_{i=1}^{m}x_{i}\bigg)^{\alpha}\prod\limits_{i=1}^{m}x_{i}^{q_{i}} dx_{1}\\
&=\int_{0}^{1} dx_{m}\cdots \int_{0}^{1-\sum\limits_{i=3}^{m} x_{i}}\prod\limits_{i=2}^{m}x_i^{q_i}\bigg(
\int_0^{1-\sum\limits_{i=2}^{m}x_i}\bigg(1-\sum\limits_{i=2}^{m}x_{i}-x_1\bigg)^{\alpha}x_1^{q_1}dx_1\bigg)dx_2\\
&=B(q_1+1, \alpha+1)\int_{0}^{1} dx_{m}\cdots \int_{0}^{1-\sum\limits_{i=3}^{m} x_{i}}\prod\limits_{i=2}^{m}x_i^{q_i}\bigg(1-\sum\limits_{i=2}^{m}x_i\bigg)^{\alpha+q_1+1}dx_2\\
&=B(q_1+1,\alpha+1)B(q_2+1,\alpha+q_1+2)\cdots B(q_m+1,\alpha+\sum\limits_{i=1}^{m-1}q_i+m-1)\\
&=\frac{\prod_{i=1}^{m}\Gamma(q_{i}+1)\Gamma(\alpha+1)}{\Gamma(\alpha+\sum_{i=1}^{m}q_{i}+m+1)}.
\end{aligned}
\end{equation*}
The proof is completed.

\begin{lem}\label{lem:2}
For any $p\in\mathbb{N}^n$, $q\in\mathbb{N}^m$ and $\alpha>-1$, then we have
\begin{equation*}
\|z^pw^q\|_{2,(-\rho)^\alpha}^2=\pi^{n+m}\frac{\prod_{i=1}^{n}\Gamma(p_i+1)\prod_{i=1}^{m}\Gamma(q_i+1)\Gamma(\alpha+1)}
{\Gamma(\alpha+m+1+|q|)[\mu(\alpha+m+|q|)]^{|p|+n}},
\end{equation*}
where $w^q$, $|q|$, $\|z^pw^q\|_{2,(-\rho)^\alpha}^2$ are given by
\begin{gather*}
w^q=\prod\limits_{i=1}^{m}w_{i}^{q_i},\quad |q|=\sum\limits_{i=1}^{m}q_i, \\
\|z^pw^q\|_{2,(-\rho)^\alpha}^2=\int_{D_{n,m}(\mu)}|z^pw^q|^2(e^{-\mu \|z\|^2}-\|w\|^2)^{\alpha} dV(z,w),
\end{gather*}
for $w=(w_1,\dots,w_m)$, $q=(q_1,\dots,q_m)$.
\end{lem}

{\bf Proof.} By definition, we have
$$\|z^pw^q\|_{2,(-\rho)^\alpha}^2=\int_{D_{n,m}(\mu)}|z^pw^q|^2(e^{-\mu \|z\|^2}-\|w\|^2)^{\alpha} dV(z,w).$$
By setting $z_j=r_je^{i\theta_j}(1\le j\le n)$, $w_l=k_le^{i\theta_l}(1\le l\le m)$, we obtain
\begin{equation*}
\|z^pw^q\|_{2,(-\rho)^\alpha}^2=(2\pi)^{n+m}\int_{
\begin{subarray}{l}
\|k\|^2<e^{-\mu\|r\|^2}\\k\ge0,r\ge0
\end{subarray}}
r^{2p+1}k^{2q+1}(e^{-\mu \sum\limits_{j=1}^{n}r_j^2}-\sum\limits_{l=1}^{m}k_l^2)^{\alpha}dr dk,
\end{equation*}
where $r=(r_1,\dots,r_n)$ ,$k=(k_1,\dots,k_m)$. By setting $s_i=r_i^2(1\le i\le n)$ and $t_j=k_j^2(1\le j\le m)$, we have
\begin{equation*}
\|z^pw^q\|_{2,(-\rho)^\alpha}^2=\pi^{n+m}\int_{
\begin{subarray}{l}
\sum\limits_{j=1}^{m}t_j<e^{-\mu\sum\limits_{i=1}^{n}s_i}\\t_j\ge0,s_i\ge0
\end{subarray}}
s^pt^q(e^{-\mu\sum\limits_{i=1}^{n}s_i}-\sum\limits_{j=1}^{m}t_j)^\alpha ds dt.
\end{equation*}
Letting $\widetilde{t_j}=e^{\mu\sum\limits_{i=1}^{n}s_i}t_j$, it follows
\begin{equation*}
\|z^pw^q\|_{2,(-\rho)^\alpha}^2=\pi^{n+m}\int_{(\mathbb{R}_+)^n}s^pe^{-\mu(\alpha+m+|q|)\sum\limits_{i=1}^{n}s_i} ds
\int_{
\begin{subarray}{l}
\sum\limits_{j=1}^{m}\widetilde{t_j}<1 \\
\widetilde{t_j}\ge0
\end{subarray}}
(1-\sum\limits_{j=1}^{m}\widetilde{t_j})^\alpha\widetilde{t}^q d\widetilde{t}.
\end{equation*}
Since $\alpha>-1$, by Lemma \ref{lem:1}, we have
\begin{equation*}
\|z^pw^q\|_{2,(-\rho)^\alpha}^2=\pi^{n+m}\frac{\prod\limits_{i=1}^{m}\Gamma(q_i+1)\Gamma(\alpha+1)}{\Gamma(\alpha+m+1+\sum\limits_{i=1}^{m}q_i)}
\int_{(\mathbb{R}_+)^n}s^p e^{-\mu(\alpha+m+|q|)\sum\limits_{i=1}^{n}s_i}ds.
\end{equation*}
Since
\begin{equation*}
\int_{0}^{\infty}s_i^{p_i}e^{-\mu(\alpha+m+|q|)s_i} ds_i=[\mu(\alpha+m+|q|)]^{-p_i-1}\Gamma(p_i+1),
\end{equation*}
we obtain
\begin{equation*}
\|z^pw^q\|_{2,(-\rho)^\alpha}^2=\pi^{n+m}\frac{\prod\limits_{i=1}^{n}\Gamma(p_i+1)\prod\limits_{i=1}^{m}\Gamma(q_i+1)\Gamma(\alpha+1)}
{\Gamma(\alpha+m+1+|q|)[\mu(\alpha+m+|q|)]^{|p|+n}}.
\end{equation*} The proof is completed.

Let $A^p(\mathbb{C}^n, \eta_\alpha)$ be the space of all entire function $f$ on $\mathbb{C}^n, n\ge 1,$ such that $|f|^p$ is integrable with respect to the Gaussian
$$\eta_\alpha(z)=e^{-\alpha\|z\|^2},$$
where $\alpha>0, 1\le p<\infty.$ Equipped with the norm inherited from $L_\alpha^p(\mathbb{C}^n,\eta_\alpha)$, $A^p(\mathbb{C}^n, \eta_\alpha)$ become Banach spaces. 
In particular, $A^2(\mathbb{C}^n, \eta_\alpha)$ is the Segal-Bargmann-Fock space of quantum mechanics with parameter $\alpha$. 
The function
\begin{equation}\label{eq:ker}
K_{\mathbb{C}^n, \eta_\alpha}(x,y)
=\bigg(\frac{\alpha}{\pi}\bigg)^ne^{\alpha\left<x,y\right>},\quad x,y\in\mathbb{C}^n,
\end{equation}
is the Bergman kernel for $A^2(\mathbb{C}^n, \eta_\alpha)$ (e.g., see \cite{Bom} for references here).

The integral operator defined by
$$P_{\mathbb{C}^n, \eta_\alpha} f(x)=\int_{\mathbb{C}^n}f(y)K_{\mathbb{C}^n, \eta_\alpha}(x,y) \eta_\alpha(y) dV(y),\quad x\in\mathbb{C}^n,$$
is the orthogonal projection in $L^2(\mathbb{C}^n, \eta_\alpha)$ onto $A^2(\mathbb{C}^n, \eta_\alpha)$. $P_{\mathbb{C}^n, \eta_\alpha}$ is bounded on $L^2(\mathbb{C}^n, \eta_\alpha)$, but this turns to be no longer the case for $L^p(\mathbb{C}^n, \eta_\alpha)$ with $p\ne 2$. Janson, Peetre and Rochberg proved the following assertions.

\begin{thm}{\rm(See \cite{Jan})}\label{thm3}
Let $\alpha\in\mathbb{R}, \beta>0, 1\le p<\infty$ satisfy $\beta p>\alpha$. Then $P_{\mathbb{C}^n, \eta_\beta }$ is bounded from $L^p(\mathbb{C}^n, \eta_\alpha)$ into $L^p(\mathbb{C}^n, \eta_\gamma)$, where $\frac{1}{\gamma}=\frac{4(\beta p-\alpha)}{p^2\beta^2}$. Inparticular, $P_{\mathbb{C}^n, \eta_\alpha}$ is bounded on $L^p(\mathbb{C}^n, \eta_\alpha)$ if and only if $ p=2$.
\end{thm}

By comparing the expression of the bergman kernel in Theorem \ref{thm2} and \eqref{eq:ker} we have the following lemma.

\begin{lem}\label{pro:1}
Let $K_{D_{n,m}(\mu),(-\rho)^\alpha}$ be the Bergman kernel for $A^2(D_{n,m}(\mu),(-\rho)^\alpha)$ and $K_{\mathbb{C}^n, \eta_\alpha}$ be the Bergman kernel for $A^2(\mathbb{C}^n, \eta_\alpha)$. Then we have
\begin{equation*}
K_{D_{n,m}(\mu),(-\rho)^\alpha}((x,0),(s,0))=\frac{\Gamma(\alpha+m+1)}{\pi^m\Gamma(\alpha+1)}K_{\mathbb{C}^n, \eta_{\mu(\alpha+m)}}(x,s),
\end{equation*}
where $x,s\in\mathbb{C}^n$.
\end{lem}

\section{Proofs of the main results}
\subsection{Proof of Theorem \ref{thm2}}
Since $\{\frac{z^pw^q}{\|z^pw^q\|_{2,(-\rho)^\alpha}}\}$ constitutes an orthonormal basis of $A^2(D_{n,m}(\mu), (-\rho)^\alpha dv)$ and the Fock-Bargmann-Hartogs domain $D_{n,m}(\mu)$ is a Reinhardt domain, we have
\begin{equation*}
K_{D_{n,m}(\mu),(-\rho)^\alpha}((x,y),(s,t))=\sum_{p\in\mathbb{N}^n,q\in\mathbb{N}^m}\frac{x^py^q\overline{s^pt^q}}{\|x^py^q\|_{2,(-\rho)^\alpha}\|s^pt^q\|_{2,(-\rho)^\alpha}}.
\end{equation*}
By Lemma \ref{lem:1}, we have

\begin{align}\label{eq:ker1}
K_{D_{n,m}(\mu),(-\rho)^\alpha}((x,y),(s,t))
&=\sum_{p\in\mathbb{N}^n,q\in\mathbb{N}^m}\frac{\Gamma(\alpha+m+1+|q|)[\mu(\alpha+m+|q|)]^{n+|p|}}
{\pi^{n+m}\prod\limits_{i=1}^{n}\Gamma(p_i+1)\prod\limits_{i=1}^{m}\Gamma(q_i+1)\Gamma(\alpha+1)}x^py^q\overline{s^pt^q}\nonumber\\
&=\frac{1}{\pi^{n+m}}\sum\limits_{q\in\mathbb{N}^m}\phi(x,s)
\frac{\Gamma(\alpha+m+1+|q|)[\mu(\alpha+m+|q|)]^n}{\prod\limits_{i=1}^{m}\Gamma(q_i+1)\Gamma(\alpha+1)}y^q\overline{t^q},
\end{align}
where $\phi(x,s)=\sum\limits_{p\in\mathbb{N}^n}\frac{[\mu(\alpha+m+|q|)]^{|p|}}{\prod\limits_{i=1}^{n}\Gamma(p_i+1)}x^p\overline{s^p}$.

It is easy to calculate that
\begin{equation}\label{eq:ker2}
\sum\limits_{p\in\mathbb{N}^n}\frac{[\mu(\alpha+m+|q|)]^{|p|}}{\prod\limits_{i=1}^{n}\Gamma(p_i+1)}x^p\overline{s^p}
=e^{\mu(\alpha+m+|q|)\left<x,s\right>}.
\end{equation}
Subsitituting \eqref{eq:ker2} into \eqref{eq:ker1}, we obtain
\begin{align*}
K_{D_{n,m}(\mu),(-\rho)^\alpha}((x,y),(s,t))
&=\frac{1}{\pi^{n+m}}\sum\limits_{q\in\mathbb{N}^m}
e^{\mu(\alpha+m+|q|)\left<x,s\right>}
\frac{\Gamma(\alpha+m+1+|q|)[\mu(\alpha+m+|q|)]^n}{\prod\limits_{i=1}^{m}\Gamma(q_i+1)\Gamma(\alpha+1)}y^q\overline{t^q}\\
&=\frac{\mu^n}{\pi^{n+m}}\sum_{q\in\mathbb{N}^m}
\frac{\Gamma(\alpha+m+1+|q|)(\alpha+m+|q|)^n}{\prod\limits_{i=1}^{m}\Gamma(q_i+1)\Gamma(\alpha+1)}
e^{\mu(\alpha+m+|q|)\left<x,s\right>}y^q\overline{t^q}\\
&=\frac{\mu^n}{\pi^{n+m}}\sum_{k\in\mathbb{N}}
\frac{\Gamma(\alpha+m+k+1)(\alpha+m+k)^n}{\Gamma(\alpha+1)\Gamma(k+1)}e^{\mu(\alpha+m+k) \left< x,s \right>}\left< y,t \right>^{k}.
\end{align*}
The proof is completed.

\subsection{Proof of Theorem \ref{th:1}}

For a given $p\in [1,\infty)\setminus\{2\}$, by Theorem \ref{thm3}, $P_{\mathbb{C}^n, \eta_{\mu(\alpha+m)}}$ is unbounded on $L^p(\mathbb{C}^n,\eta_{\mu(\alpha+m)})$, where $\eta_{\mu(\alpha+m)}=e^{-\mu(\alpha+m)\|z\|^2}$.
Therefore, there exists a sequence $\{f_n(z)\}$ in $L^p(\mathbb{C}^n,\eta_{\mu(\alpha+m)})
\cap L^2(\mathbb{C}^n,\eta_{\mu(\alpha+m)})$  such that
\begin{equation}\label{eq:001}
\lim_{n\to\infty}\frac{\|P_{\mathbb{C}^n,\eta_{\mu(\alpha+m)}} f_n\|_{p,\eta_{\mu(\alpha+m)}}^p}{\|f_n\|_{p,\eta_{\mu(\alpha+m)}}^p}=\infty.
\end{equation}

Define $F_n(z,w)=f_n(z)$. Then
\begin{equation*}
\begin{aligned}
\|F_n\|_{p,(-\rho)^\alpha}^p
&=\int_{D_{n,m}(\mu)}|F_n(z,w)|^p(e^{-\mu\|z\|^2}-\|w\|^2)^\alpha dV(z,w)\\
&=\int_{\mathbb{C}^n}|f_n(z)|^pe^{-\mu\alpha\|z\|^2}\bigg(\int_{\|w\|^2<e^{-\mu\|z\|^2}}(1-e^{\mu\|z\|^2}\|w\|^2)^\alpha dV(w)\bigg) dV(z).
\end{aligned}
\end{equation*}
Let $\sigma$ be the rotation-invariant positive Borel measure on $\partial \mathbb{B}^m$, the surface of unit ball of complex dimension m, with $\sigma(\partial \mathbb{B}^m)=1$, and let $w=r\zeta,\; \zeta\in\partial\mathbb{B}^m$. Then we have
\begin{equation*}
\begin{aligned}
\|F_n\|_{p,(-\rho)^\alpha}^p
&=\int_{\mathbb{C}^n}|f_n(z)|^pe^{-\mu\alpha\|z\|^2}\int_{0}^{e^{-\frac{\mu}{2}\|z\|^2}}2mV(\mathbb{B}^m)r^{2m-1}(1-e^{\mu\|z\|^2}r^2)^\alpha dr dV(z)\\
&=mB(m,\alpha+1)V(\mathbb{B}^m)\int_{\mathbb{C}^n}|f_n(z)|^pe^{-\mu(\alpha+m)\|z\|^2} dV(z)\\
&=mB(m,\alpha+1)V(\mathbb{B}^m)\|f_n\|_{p,\eta_{\mu(\alpha+m)}}^p,
\end{aligned}
\end{equation*}
where $B(m,\alpha+1)$ is the beta function. Therefore, $F_n(z,w)\in L^p(D_{n,m}(\mu),(-\rho)^\alpha)$ for any $n$, and
\begin{equation}\label{eq:1}
\|F_n\|_{p,(-\rho)^\alpha}^p=mB(m,\alpha+1)V(\mathbb{B}^m)\|f_n\|_{p,\eta_{\mu(\alpha+m)}}^p.
\end{equation}
\begin{align}
&P_{D_{n,m}(\mu),(-\rho)^\alpha}F_n(z,0)\nonumber\\
&=\int_{D_{n,m}(\mu)}K_{D_{n,m}(\mu),(-\rho)^\alpha}((z,0),(s,t))F_n(s,t)(-\rho(s,t))^\alpha dV(s,t)\nonumber\\
&=\int_{\mathbb{C}^n}\int_{\|t\|^2<e^{-\frac{\mu}{2}\|s\|^2}}K_{D_{n,m}(\mu),(-\rho)^\alpha}((z,0),(s,t))f_n(s)(-\rho(s,t))^\alpha dV(t) dV(s)\nonumber\\
&=\int_{\mathbb{C}^n}\bigg(\int_{0}^{e^{-\frac{\mu}{2}\|s\|^2}}2mV(\mathbb{B}^m)r^{2m-1} dr \nonumber\\
&\quad \times  \int_{\partial\mathbb{B}^m}K_{D_{n,m}(\mu),(-\rho)^\alpha}((z,0),(s,r\zeta))f_n(s)(e^{-\mu\|s\|^2}-r^2)^\alpha d\sigma(\zeta)\bigg) dV(s)\nonumber\\
&=\int_{\mathbb{C}^n}\bigg(2mV(\mathbb{B}^m)f_n(s)\int_{0}^{e^{-\frac{\mu}{2}\|s\|^2}}r^{2m-1}(e^{-\mu\|s\|^2}-r^2)^\alpha dr \nonumber\\
&\quad \times
\int_{\partial\mathbb{B}^m}K_{D_{n,m}(\mu),(-\rho)^\alpha}((z,0),(s,r\zeta)) d\sigma(\zeta)\bigg) dV(s). \label{eq:2}
\end{align}
Since $K_{D_{n,m}(\mu),(-\rho)^\alpha}((z,w),(s,t))$ is antiholomorphic in $t$, by the mean value property we have
\begin{equation}\label{eq:3}
\int_{\partial\mathbb{B}^m}K_{D_{n,m}(\mu),(-\rho)^\alpha}((z,0),(s,r\zeta)) d\sigma(\zeta)=V(\partial\mathbb{B}^m)K_{D_{n,m}(\mu),(-\rho)^\alpha}((z,0),(s,0)),
\end{equation}
where $V(\partial\mathbb{B}^m)$ is the volume of $\partial\mathbb{B}^m$.
Putting \eqref{eq:3} into \eqref{eq:2}, we obtain
\begin{align}
&P_{D_{n,m}(\mu),(-\rho)^\alpha}F_n(z,0)\nonumber\\
&=\int_{\mathbb{C}^n}\bigg(2mV(\partial\mathbb{B}^m)V(\mathbb{B}^m)K_{D_{n,m}(\mu),(-\rho)^\alpha}((z,0),(s,0))f_n(s)\int_{0}^{e^{-\frac{\mu}{2}\|s\|^2}}
r^{2m-1}(e^{-\mu\|s\|^2}-r^2)^\alpha dr\bigg) dV(s)\nonumber \\
&=mV(\partial\mathbb{B}^m)V(\mathbb{B}^m)B(m,\alpha+1)\int_{\mathbb{C}^n}K_{D_{n,m}(\mu),(-\rho)^\alpha}((z,0),(s,0))f_n(s)e^{-\mu(\alpha+m)\|s\|^2} dV(s).\label{eq:4}
\end{align}
Applying Lemma \ref{pro:1} to \eqref{eq:4}, we get
\begin{align}
P_{D_{n,m}(\mu),(-\rho)^\alpha}F_n(z,0)
&=c
\int_{\mathbb{C}^n}K_{\mathbb{C}^n, \eta_{\mu(\alpha+m)}} (z,s)f_n(s)e^{-\mu(\alpha+m)\|s\|^2} dV(s)\nonumber\\
&=cP_{\mathbb{C}^n, \eta_{\mu(\alpha+m)}} f_n(z),\label{eq:4.1}
\end{align}
where $c=\frac{mV(\partial\mathbb{B}^m)V(\mathbb{B}^m)B(m,\alpha+1)\Gamma(\alpha+m+1)}{\pi^m\Gamma(\alpha+1)}$.

Next we estimate the norm of $P_{D_{n,m}(\mu),(-\rho)^\alpha}F_n$.
\begin{align}
&\|P_{D_{n,m}(\mu),(-\rho)^\alpha}F_n\|_{p,(-\rho)^\alpha}^p\nonumber\\
&=\int_{D_{n,m}(\mu)}|P_{D_{n,m}(\mu),(-\rho)^\alpha}F_n(z,w)|^p(e^{-\mu\|z\|^2}-\|w\|^2)^\alpha dV(z,w)\nonumber\\
&=\int_{\mathbb{C}^n}\int_{\|w\|^2<e^{-\mu\|z\|^2}}|P_{D_{n,m}(\mu),(-\rho)^\alpha}F_n(z,w)|^p(e^{-\mu\|z\|^2}-\|w\|^2)^\alpha dV(w) dV(z)\nonumber\\
&=\int_{\mathbb{C}^n}\bigg(\int_{0}^{e^{-\frac{\mu}{2}\|z\|^2}}2mV(\mathbb{B}^m)r^{2m-1}(e^{-\mu\|z\|^2}-r^2)^\alpha dr \int_{\partial\mathbb{B}^m}|P_{D_{n,m}(\mu),(-\rho)^\alpha}F_n(z,r\zeta)|^p d\sigma(\zeta)\bigg)dV(z). \label{eq:5}
\end{align}
Since $P_{D_{n,m}(\mu),(-\rho)^\alpha}F_n(z,w)$ is holomorphic in $w$, by the submean value property we have
\begin{equation}\label{eq:6}
\int_{\partial\mathbb{B}^m}|P_{D_{n,m}(\mu),(-\rho)^\alpha}F_n(z,r\zeta)|^p d\sigma(\zeta)\ge
V(\partial\mathbb{B}^m)|P_{D_{n,m}(\mu),(-\rho)^\alpha}F_n(z,0)|^p.
\end{equation}
Subsitituting \eqref{eq:6} in \eqref{eq:5}, together with \eqref{eq:4.1},  we obatin
\begin{align}
&\|P_{D_{n,m}(\mu),(-\rho)^\alpha}F_n\|_{p,(-\rho)^\alpha}^p\nonumber\\
&\ge\int_{\mathbb{C}^n}\bigg(\int_{0}^{e^{-\frac{\mu}{2}\|z\|^2}}2mV(\mathbb{B}^m)r^{2m-1}(e^{-\mu\|z\|^2}-r^2)^\alpha V(\partial\mathbb{B}^m)|P_{D_{n,m}(\mu),(-\rho)^\alpha}F_n(z,0)|^p dr\bigg) dV(z)\nonumber\\
&=mB(m,\alpha+1)V(\partial\mathbb{B}^m)V(\mathbb{B}^m)c^p\int_{\mathbb{C}^n}|P_{\mathbb{C}^n,\eta_{\mu(\alpha+m)}} f_n(z)|^pe^{-\mu(\alpha+m)\|z\|^2}dV(z)\nonumber\\
&=mB(m,\alpha+1)V(\partial\mathbb{B}^m)V(\mathbb{B}^m)c^p\|P_{\mathbb{C}^n,\eta_{\mu(\alpha+m)}}f_n\|_{p,\eta_{\mu(\alpha+m)}}^p.\label{eq:7}
\end{align}

By \eqref{eq:1} and \eqref{eq:7} we know that
\begin{equation}\label{eq:007}
\frac{\|P_{D_{n,m}(\mu),(-\rho)^\alpha}F_n\|_{p,(-\rho)^\alpha}^p}{\|F_n\|_{p,(-\rho)^\alpha}^p}\ge V(\partial\mathbb{B}^m)c^p
\frac{\|P_{\mathbb{C}^n,\eta_{\mu(\alpha+m)}} f_n\|_{p,\eta_{\mu(\alpha+m)}}^p}{\|f_n\|_{p,\eta_{\mu(\alpha+m)}}^p}.
\end{equation}
Thus, by \eqref{eq:001} and \eqref{eq:007}, we get that $\lim\limits_{n\to\infty}\frac{\|P_{D_{n,m}(\mu),(-\rho)^\alpha}F_n\|_{p,(-\rho)^\alpha}^p}{\|F_n\|_{p,(-\rho)^\alpha}^p}= \infty$.
This means that $P_{D_{n,m}(\mu),(-\rho)^\alpha}$ is unbounded on $L^p(D_{n,m}(\mu), (-\rho)^\alpha)$ for $p\in [1,\infty)\setminus\{2\}$.

Therefore, $P_{D_{n,m}(\mu),(-\rho)^\alpha}$ is bounded on $L^p(D_{n,m}(\mu), (-\rho)^\alpha)$ if and only if $p=2$.
The proof is completed.

\vskip 10pt

\noindent\textbf{Acknowledgments}\quad  The project is supported by the National Natural Science Foundation of China (No. 11671306).

\addcontentsline{toc}{section}{References}
\phantomsection
\renewcommand\refname{References}

\clearpage
\end{document}